\documentclass[a4paper,12pt]{amsart}

\nonstopmode

\newcommand{\version}{10}

\title[D-minimal groups v.~\version]{Groups and rings definable in d-minimal structures\\
\small{Version~\version}}
\author[A. Fornasiero]{Antongiulio Fornasiero}
\address{Institut f\"ur Mathematische Logik\\
    Einsteinstr.~62, 48149 M\"unster, Germany}
\email{antongiulio.fornasiero@googlemail.com} 

\newcommand{\mykeywords}{D-minimal, o-minimal, group, locally o-minimal, definably complete}

\date{22 Sep 2012}

\usepackage{ifpdf}

\usepackage{amsmath, amsthm, amssymb}
\usepackage{mathtools}
\usepackage{xspace}

\usepackage{microtype}

\usepackage{comment}

\usepackage[neveradjust]{paralist}
\AtBeginDocument{\setdefaultleftmargin{0pt}{}{}{}{}{}}
\setdefaultenum{\upshape(1)}{\upshape i)}{}{}

\ifpdf
\usepackage[pdftex,
pdfborder={0 0 0}, plainpages=false,%
pdfpagelabels, pdfdisplaydoctitle=true,%
pdfstartview={XYZ null null null},%
colorlinks=false
]{hyperref}
\hypersetup{%
  pdfauthor={A. Fornasiero},
  pdftitle={D-minimal groups},
  pdfkeywords={\mykeywords},
  pdfsubject={Groups definable in d-minimal structures}
  }
\else
\usepackage[plainpages=false,linktocpage=true,pdfpagelabels, colorlinks=false
]{hyperref}
\fi

\usepackage[backrefs, alphabetic, msc-links]{amsrefs}

\makeatletter
\@namedef{subjclassname@2010}{%
 \textup{2010} Mathematics Subject Classification}
\makeatother

\newcommand*{\intro}[1]{\textbf{#1}}
\newcommand*{\Pa}[1]{\bigl( #1 \bigr)}
\newcommand*{\set}[1]{\{\,#1\,\}}
\newcommand*{\mset}[1]{\{#1\}}

\newcommand{\N}{\mathbb{N}}
\newcommand{\Real}{\mathbb{R}}
\newcommand{\R}{\mathbb{R}}
\newcommand{\Raz}{\mathbb{Q}}
\newcommand{\Qs}{\mathbb{Q}^*}

\newcommand{\Rbar}{\bar\Real}
\DeclareMathOperator{\GL}{GL}
\DeclareMathOperator{\SL}{SL}

\newcommand*{\cl}[1]{\overline{#1}}
\DeclareMathOperator{\cll}{cl}

\newcommand{\rest}{\upharpoonright}

\DeclareMathOperator{\id}{id}

\newcommand*{\pair}[1]{\langle #1 \rangle}
\newcommand{\Am}{\mathbb A}

\newcommand{\Dm}{\mathbb D}
\newcommand{\Lm}{\mathbb L}

\newcommand{\av}{{\bar a}}

\newcommand{\cv}{{\bar c}}

\newcommand{\x}{{\bar x}}
\newcommand{\y}{{\bar y}}

\def\Ind#1#2{#1\setbox0=\hbox{$#1x$}\kern\wd0\hbox to
  0pt{\hss$#1\mid$\hss}\lower.9\ht0\hbox to 0pt{\hss$#1\smile$\hss}\kern\wd0}

\newcommand{\rplus}{\oplus}

\newcommand{\K}{\mathbb{K}}
\newcommand{\Z}{\mathbb{Z}}
\newcommand{\Rc}{\mathcal R}
\newcommand{\F}{\mathbb{F}}
\newcommand{\G}{\mathbb{G}}
\newcommand{\Hg}{\mathbb{H}}

\newcommand{\Kb}{\overline{\K}}

\newcommand{\Cp}{\mathcal C^p}
\newcommand{\Ck}{\mathcal C^k}
\newcommand{\Cone}{\mathcal C^1}
\newcommand{\Cinf}{\mathcal C^\infty}

\DeclareMathOperator{\Ker}{Ker}

\DeclareMathOperator{\Ad}{Ad}
\DeclareMathOperator{\ad}{ad}
\newcommand{\lZ}{\mathfrak Z}
\newcommand{\hl}{\mathfrak h}
\newcommand{\gl}{\mathfrak g}
\newcommand{\Gz}{\G^{0}}
\newcommand{\Gzz}{\G^{00}}

\DeclareMathOperator{\reg}{reg}

\newcommand{\ominimal}{o\hyph minimal\xspace}

\newcommand{\dminimal}{d\hyph minimal\xspace}

\newcommand{\zdvs}{zero\hyph divisors\xspace}
\newcommand{\nzdv}{nonzero\hyph divisor\xspace}

\DeclareMathOperator{\RK}{rk}
\newcommand{\zrk}{\RK^Z}

\def\hyph{\nobreakdash-\hspace{0pt}\relax}
\providecommand{\rom}{\textup}

\newcommand{\Wlog}{W.l.o.g\mbox{.}\xspace}
\newcommand{\wloG}{w.l.o.g\mbox{.}\xspace}
\newcommand{\eg}{e.g\mbox{.}\xspace}
\newcommand{\ie}{i.e\mbox{.}\xspace}

\newcommand{\wrt}{w.r.t\mbox{.}\xspace}

\newcommand{\tfae}{t.f.a.e\mbox{.}\xspace}
\newcommand{\Tfae}{T.f.a.e\mbox{.}\xspace}

\newcommand{\resp}{resp\mbox{.}\xspace}

\newtheorem{lemma}{Lemma}[section]
\newtheorem{thm}[lemma]{Theorem}
\newtheorem{corollary}[lemma]{Corollary}
\newtheorem{conjecture}[lemma]{Conjecture}
\newtheorem{proposition}[lemma]{Proposition}
\newtheorem{open problem}[lemma]{Open problem}

\newtheorem*{fact*}{Fact}
\newtheorem{fact}[lemma]{Fact}

\theoremstyle{remark}

\newtheorem{claim}{Claim}
\newtheorem*{claim*}{Claim}

\theoremstyle{definition}
\newtheorem{definizione}[lemma]{Definition}

\newtheorem*{notation*}{Notation}
\newtheorem{remark}[lemma]{Remark}
\newtheorem{final remark}[lemma]{Final remark}
\newtheorem{example}[lemma]{Example}
\newtheorem{examples}[lemma]{Examples}
\newtheorem{question}[lemma]{Question}

\newtheorem{warning}[lemma]{Warning}

\begin{document}

\begin{abstract}
We study groups and rings definable in d-minimal expansions of ordered fields.
We generalize to such objects some known results from o-minimality.
In particular, we prove that we can endow a definable group with a unique
definable topology making it a definable manifold and a topological group, 
and that a definable ring of dimension at least 1 and without zero divisors is a
skew field.
\end{abstract}

\maketitle

{\small
\setcounter{tocdepth}{1}
\microtypesetup{protrusion=false}
\tableofcontents
\microtypesetup{protrusion=true}
}

\makeatletter
\renewcommand\@makefnmark%
   {\normalfont(\@textsuperscript{\normalfont\@thefnmark})}
\renewcommand\@makefntext[1]%
   {\noindent\makebox[1.8em][r]{\@makefnmark\ }#1}
\makeatother

\section{Introduction}

Let $\K$ be a first-order expansion of an ordered field.
Recall that $\K$ is \ominimal if every definable
(by ``definable'' we will always mean ``definable in $\K$ with parameters'')
subset of $\K$ is a finite union of points an intervals with endpoints in $\K
\cup \set{\pm \infty}$.
O-minimal structures have been widely studied (see \cite{vdd} for an
introduction). 
There is a rich literature on groups definable in o-minimal structures; see
for instance~\cites{pillay88, OPP, PPS00} for problems treated in this article 
(we will discuss those references more in details in the following sections),  
and \cites{otero08, peterzil-survey10} for an overview.
One of the starting points was A.~Pillay's theorem that any such group can be
endowed in a unique way with a topology that makes it both a topological group
and a definable manifold.
An important result in this context is the solution of the ``Pillay's conjecture''  and its
extension to groups definable in NIP structures (see \cite{HPP08}).

In this article, we extend some of the previous results about groups and rings
definable in \ominimal structures to groups and rings definable in \dminimal structures.

Recall that $\K$ is said to be \intro{definably complete} if every definable
subset of $\K$ has a supremum in $\K \cup \set {\pm \infty}$
(see \eg \cite{FS:1} and its bibliography).
In this article we study the following generalization of o-minimality:

\begin{definizione}\label{def:dmin-intro}
$\K$ is \intro{\dminimal{}} if it is definably complete, and every definable
set $X \subset \K$ is the union of an open set and finitely many discrete sets,
where the number of discrete sets does not depend on the parameters of
definition of~$X$.
\end{definizione}
Notice that $\K$ is \ominimal iff it is definably complete and every definable
unary set is the union of an open set and finitely many points.

\cite{dries85} gives the first known example of a \dminimal non o-minimal
expansion of~$\Real$, 
\cites{FM05, MT06} give more examples of d-minimal expansions of~$\Real$ (and
introduce the notion of d-minimality), 
\cite{miller05} studies general properties of \dminimal expansions of~$\Real$
(and other such ``tameness'' notions),
and in \cite{fornasiero:d-min} we studied \dminimal structures in general
(see \S\ref{sec:preliminary} for an overview).

Let $\K$ be a \dminimal structure.
The main results are:
\begin{thm}[\cite{pillay88}]
Let $\G \coloneqq \pair{G, \cdot, e, {}^{-1}}$ be a definable group.
Then, we can endow $G$ with a unique differential structure 
\rom(see Definition~\ref{def:manifold}\rom), such that $\cdot$ and $^{-1}$ are
differentiable functions.

Similarly, if $\F$ is a definable ring, then  
we can endow $F$ with a unique differential structure, such that the ring
operations are differentiable functions. 
\end{thm}

\begin{thm}
Let $\G$ be a definable Abelian group.
If $\G$ is definably connected, then it is divisible.
\end{thm}

\begin{thm}[\cite{PPS00}]
Let $\G$ be a definable group.
Assume that $\G$ is definably connected, centerless, and semisimple 
\rom(\ie, every definable normal Abelian subgroup is discrete\rom).
Then, $\G$ is definably isomorphic to a semi-algebraic group.
\end{thm}


\begin{thm}[\cite{OPP}]
Let $\F$ be a definable ring without 0 divisors, such that $\dim(\F) \geq 1$.
Then, $\F$ is a skew field, and it is definably isomorphic to either $\K$,
$\K(\sqrt{-1})$,  or the rings of quaternions over~$\K$.
\end{thm}

As far as we know, the following theorem is new even for o-minimal structures.

\begin{thm}
Let $\F$ be a definable ring with~$1$.
Then, there exist a definable 
$\K$-algebra $\F^0 \subseteq \F$ and a definable discrete
subring $D \subseteq \F$, such that, as rings,
$\F = D \rplus \F^0$.
\end{thm}


\subsection*{Acknowledgments}
Thanks to A.~Berarducci, L.~Kramer, and A.~Pillay  for the useful discussions .


\section{Preliminaries on d-minimal structures}
\label{sec:preliminary}

\subsection{Conventions and notation}
$\K$ will always be a \dminimal expansion of an ordered field
$\Kb \coloneqq \pair{K, +, \cdot, <, 0, 1}$, and
``definable'' will always mean ``definable with parameters from~$K$''.
Moreover, $\cl X$ or $\cll(X)$ denote the topological closure of~$X$.
$\Rbar \coloneqq \pair{\R, +, \cdot, <, 0, 1}$ is the field of real numbers.


\subsection{Some examples}
\label{sec:dmin-example}
\begin{enumerate}
\item A definably complete expansion of an ordered field $\K$ 
is \intro{locally o-minimal}
if every definable subset of $K$ with empty interior is discrete
(see \cites{fornasiero-lomin, schoutens}).
Clearly, a locally \ominimal structure is \dminimal, and
an ultraproduct of \ominimal structures is locally \ominimal (but not
necessarily \ominimal).
\item Given $c \in \R$, denote $c^\Z \coloneqq \set{c^n: n \in \Z}$.
The structure $\pair{\R, +, \cdot, c^\Z}$ is \dminimal for every $c \in \R$
(see \cite{dries85});
on the other, the structure $\pair{\R, +, \cdot, 2^\Z, 3^\Z}$ not only is not
\dminimal, but it defines the set of natural numbers
(see \cite{hieronymi}).
\item\label{ex:fast} Let $\mathcal R$ be an \ominimal expansion of $\Rbar$, 
$(a_i)_{i \in \N}$ be a ``fast sequence'',
let $A \coloneqq \set{a_i: i \in \N}$.
Then,
the expansion of $\Rbar$ by a predicate for each subset
$A^n$ for every $n \in \N$ is \dminimal
(see \cite{FM05} for the relevant definitions and proofs).
\item For other examples of \dminimal expansions of $\Rbar$, see
\cite{miller05}.
\end{enumerate}

Notice that NIP can fail in \dminimal structures: there are \dminimal
expansion of $\Rbar$ that define an isomorphic copy of $\pair{\N, +, \cdot}$
(see Example~\ref{ex:fast})
and there exist ultraproducts of \ominimal structures that do not satisfy NIP
(see \cite{fornasiero-lomin}).


\subsection{Previous results}

See \cite{miller05} for \dminimal expansions of~$\Rbar$,
and \cite{fornasiero:d-min} for general \dminimal structures.

\begin{definizione}
Let $X \subseteq \K^n$ be a definable set.
We say that $X$ is an \intro{embedded $\Ck$-manifold} of dimension~$d$ if, for every
$x \in X$ there exists an open box $U$ containing $x$,
such that, after a permutation of coordinates,
$X \cap U$ is the graph of a $\Ck$-function $f: V \to W$, where $V \subset
\K^d$ and $W \subset \K^{n - d}$ are the unique boxes satisfying 
$U = V \times W$. 
\end{definizione}

Here we will use the following facts (from \cite{fornasiero:d-min}).

\begin{fact}
Let $X$ be a definable set.
Then, $X$ is the disjoint union of finitely many embedded manifolds.
\end{fact}

\begin{definizione}
Let $X \subseteq \K^n$ be a definable set.
The \intro{dimension} of $X$, denoted by $\dim(X)$, is $- \infty$ if $X$ is
empty, and otherwise it
is the smallest integer $d$ such
that there exists a $d$-dimensional coordinate space~$L$,
such that $\pi_L(X)$ has nonempty interior in~$L$,
where $\pi(L)$ is the orthogonal projection onto~$L$.
\end{definizione}

\begin{fact}
$\dim$ satisfies the axioms for a dimension function in~\cite{dries89}.
In particular, $\dim(X \cup Y) = \max(\dim(X), \dim(Y))$, and
if $f: X \to Y$ is a definable function, such that
$\dim(f^{-1}(y)) = d$ for every $y \in Y$, then $\dim(Y) = d + \dim(X)$.
Moreover, if $X$ is an embedded manifold, then $\dim(X)$ coincides with the
dimension of $X$ as a manifold.
Finally, $\dim(X) \leq 0$ iff $X$ is a finite union of \rom(definable\rom)
discrete sets.
\end{fact}
However, unlike the \ominimal case,
it is not true in general that, for $X$ nonempty,
$\dim(\cl X \setminus X) < \dim(X)$.

Notice that every embedded manifold is locally closed in $\K^n$.
Thus, we have the following fact.
\begin{fact}
Let $X \subseteq \K^n$ be definable.
Then, $X$ is constructible: that is, $X$ is the union of finitely many
definable locally closed sets.
\end{fact}

\begin{fact}
Let $f: \K^n \to \K^m$ be a definable function, and $k \in \N$.
Then, there exists a closed definable set $C \subset \K^n$ with empty
interior, such that, outside $C$, $f$ is~$\Ck$.
\end{fact}

For every $\Cone$ function $f : \K^m \to \K^n$, define
$\Lambda_f(k) \coloneqq \set{x \in \K^m: \RK \Pa{d_x(f))} \leq k}$, 
(where $\RK(M)$ is the rank of the matrix~$M$ and $d_x(f)$ is the differential
of $f$ at~$x$), and $\Sigma_f(k) \coloneqq f\Pa{\Lambda_f(k)}$.
The set of singular values of $f$ is 
$\Sigma_f \coloneqq \bigcup_{k = 0}^{n - 1} \Sigma_f(k)$.

\begin{fact}[Sard's Lemma]\label{fact:Sard} 
Let $f : \K^m \to \K^n$ be definable and~$\Cone$.
Then, for every $d \leq n$, $\dim\Pa{\Sigma_f(d)} \leq d$.
In particular, there exists $c \in \K^m$ that is a regular value for~$f$.
\end{fact}

\begin{definizione}
Let $X \subseteq \K^n$ be a definable set, and $p \in \N$.
Let $\reg^p(X)$ denote the set of all $x \in X$ such that, for 
some open box~$U$ containing~$a$, $X \cap U$ is a $\Cp$ embedded
manifold of the same dimension as~$X$.
\end{definizione}

\begin{fact}\label{fact:i-regular-2}
Let $X \subseteq \K^n$ be a definable set, and $p \in \N$.
Then, $X \setminus \reg^p(X)$ is nowhere dense in~$X$.
\end{fact}

\begin{fact}[Dimension is local: see \cite{FH})]
Let $X \subseteq \K^n$ be a definable set.
Assume that, for every $x \in X$, there exists a definable neighbourhood $U$
of $x$, such that $\dim(U \cap X) \leq d$.
Then, $\dim(X) \leq d$.
\end{fact}

\begin{fact}
$\K$ has definable Skolem functions and definable choice.
\end{fact}
Thus, (almost) all the results about definable sets can be extended to sets
that are interpretable in~$\K$.

\begin{definizione}
\begin{enumerate}
\item 
Let $\cv \in \K^n$ and $A \subseteq \K$.
We say that $\zrk(\cv/A) \leq d$ if there exists a $d$-dimensional set
$X$ definable with parameters from $A$, such that $\cv \in X$.
We say that $\zrk(\cv/A) = d$ if $\zrk(\cv/A) \leq d$ and 
$\zrk(\cv/A) \nleq d-1$.
\item
Let $X \subseteq \K^n$ be a definable set, 
$A \subseteq \K$ be a set containing the
parameters of definition of~$X$, and $\cv \in X$.  
We say that $\cv$ is \intro{generic} in $X$ over~$A$, 
if for every set $Y$ definable
over~$A$, if $\cv \in Y$, then $\dim(Y) \geq \dim(X)$.
\end{enumerate}
\end{definizione}

\begin{fact}
\begin{enumerate}
\item 
$\cv$ is generic in $X$ over $A$ iff $\zrk(\cv/A) = \dim(X)$.
\item If $\K$ is sufficiently saturated and $X$ is an embedded manifold, then
the set of points in $X$ that are generic over $A$ is \rom(topologically\rom)
dense in~$X$ \rom(notice that, even for o-minimal structures, this may not be
true when $X$ is not an embedded manifold\rom).
\item
$\zrk$ is the rank corresponding to a \rom(unique\rom) matroid on $\K$
\rom(see \cite{fornasiero-matroids}\rom).
\end{enumerate}
\end{fact}
Thus, most of the proofs in \ominimal situations that rely on generic elements
can be transferred without much difficulty to \dminimal structures.%
\footnote{For definable groups, there is a different notion of ``generic''
  (see~\cite{HPP08}), but we will not use it.}

\subsection{Functions}

The results in this subsection (with the same proofs) hold not only for 
\dminimal structures,
but also when $\K$ is any definably complete expansion of some ordered field.

\begin{proposition}\label{prop:differential-equation}
Let $U \coloneqq I_1 \times I_2 \subseteq \K^n \times \K^m$ be an open rectangular box,
and $F: U \to M_n(\K)$ be a $\Cone$ definable function \rom(where $M_n(\K)$ is the
set of $n \times n$ matrices over~$\K$\rom).
For $\pair{a, b} \in U$, consider the system of differential equations
\begin{equation}\label{eq:differential}
\begin{aligned}
\phi(a) &= b,\\
d_x(\phi) &= F(x, \phi(x)).
\end{aligned}
\end{equation}
Then, there exists \emph{at most} one function $\phi: I_1 \to I_2$, which is
definable, $\Cone$, and satisfies 
\eqref{eq:differential}.
\end{proposition}
\begin{proof}
The same as~\cite{OPP}*{Theorem~2.3}.
\end{proof}

\begin{proposition}[Local Submersion Theorem]
Let $f: \K^n \to \K^m$ be a definable function.
Assume that $d_0(f)$ has rank $m$.
Then, $f$ is an open map in a neighbourhood of~$0$.
\end{proposition}
\begin{proof}
The Implicit Function Theorem for definable functions 
was proved \eg in \cite{servi:PhD}.
The Local Submersion Theorem follows in the usual way.
\end{proof}


\section{Definable groups}

\begin{notation*}
In all the article, unless explicitly said otherwise, when we say that $\G$ is
a group, we will denote by $G$ the underlying set, by $\cdot$ the
multiplication, by $e$ the identity, and by ${}^{-1}$ the inverse operation
of~$\G$.
\end{notation*}

\subsection{Examples}\label{subsec:group-example}
1) Let $\G$ be a semi-algebraic group (\ie, a real Lie group definable in
$\Rbar$); let $\tilde \G \to \G$ be the universal cover of~$\G$, and let $D$ be
the kernel of the covering map $\tilde \G \to \G$.  Notice that $\Dm :=
\pair{D, +}$ is a discrete central subgroup of~$\tilde\G$.
Remember that the structure $\Rc \coloneqq \pair{\Rbar, 2^\Z}$ is \dminimal.
Since $\Dm$ is a finitely generated Abelian group, it is definable in $\Rc$.
By \cite{HPP11}*{\S8.1}, an isomorphic copy of the group $\tilde \G$, together
with the extension maps $\Dm \to \tilde \G \to \G$, can be 
interpreted in the structure 2-sorted structure $\pair{\Rbar, \Dm}$, 
and therefore the extension $\Dm \to \tilde \G \to \G$ can be interpreted in~$\Rc$.
For instance, an isomorphic copy of the group $\widetilde{\SL_2(\R)}$, together
with the extension $\Z \to \widetilde{\SL_2(\R)} \to \SL_2(\R)$,
can be defined inside $\Rc$.

2)
What can we say about groups definable in a \dminimal structure?
For $0$-dimensional groups, almost nothing.
Let $\G$ be a countable abstract group;
then, there exists a \dminimal expansion of $\Rbar$ that defines a
$0$-dimensional group that is isomorphic to~$\G$
(see \S\ref{sec:dmin-example}\eqref{ex:fast}).

The situation changes drastically for higher dimensional groups: more precisely, we will show that a group definable in a \dminimal expansion of $\Rbar$ is, in a canonical way, a Lie group.

\subsection{Definable topology}

A.~Pillay (\cite{pillay88}) proved that every group definable in an o-minimal
structure can be endowed with a definable topology that makes it both a
topological group and a definable manifold.
His result was later extended by A.~Mosley and R.~Wencel to other kinds of
structures.
While those generalizations do not apply directly to \dminimal structures, a
small modification will suffice for our purpose.

\begin{definizione}\label{def:manifold}
A \intro{finitary $\Ck$-manifold} (or, simply, finitary manifold when 
\mbox{$k = 0$}) of dimension $d$ is given by 
\begin{enumerate}\item 
a definable set $X \subseteq \K^n$, 
\item a topology $\tau$ on $X$, 
\item
finitely many definable embedded $\Ck$-manifolds $C_1, \dotsc, C_m$ of
dimension~$d$,
\item for each $i = 1, \dotsc, m$, a definable map $f_i: C_i \to X$, such
that:
\begin{enumerate}
\item 
for each $i$, $f_i$ is a homeomorphism with its image (with topology~$\tau$) 
and the image of $f_i$ is $\tau$-open in~$X$;
\item 
$X = \bigcup_i f(C_i)$
\item
for each $i$ and $j$, the partial function $f_i^{-1} \circ f_j: C_j \to C_i$
is $\Ck$ (notice that this condition is superfluous when $k = 0$).
\end{enumerate}
\end{enumerate}
We also say that $\pair{X, \tau}$ is a finitary $\Ck$-manifold if we can
find $C_1, \dotsc, C_m$ and $f_1, \dotsc, f_k$ as above;
the set $\set{\pair{f(C_1), f_1^{-1}}, \dotsc, \pair{f(C_m), f_m^{-1}}}$ is an
atlas for the finitary manifold, and a \intro{$\Ck$ differential structure}.
Two $\Ck$ differential structures on the same set $X$ are to be considered
equal if the identity map on $X$ is a $\Ck$ diffeomorphism from the first
differential structure to the second. 
We will often denote in the same way  (say, $\tau$) 
the $\Ck$ differential structure and the underlying topology.
\end{definizione}
Notice that we are not taking position on what is the ``correct'' definition
of an abstract definable manifold (\ie, one that may require infinitely many charts).
Notice also that, unlike the o-minimal case, we are \emph{not} claiming that
we can take the set $C_i$ to be open subsets of~$\K^d$; and, in fact, in
general we cannot do that: for instance, if $X$ is an infinite definable
discrete subset of~$\K$, then (obviously) 
$X$ is not a finite union of copies of~$\K^0$.

\begin{examples}
\begin{enumerate}
\item Every definable embedded manifold is a finitary manifold.
\item The disjoint union of finitely many finitary manifold of the
same dimension (with
the differential structure of  disjoint union) is also a finitary manifold.
\item The main result of \cite{fornasiero-enumerable} can be formulated by
saying that a definably complete structure with the discrete topology is
\emph{not} a finitary manifold.
It is the ``definable counterpart'' of the fact that $\R$ with the discrete
topology is not a manifold (since it is not second countable).
\end{enumerate}
\end{examples}

\begin{definizione}
Fix $k \geq 0$.
Let $X\subseteq \K^n$  be a definable set, endowed with a $\Ck$ differential
structure. 
A subset $Y \subseteq X$ is \intro{large} in $X$ if 
$\dim(X \setminus Y) < \dim X$.
We say that $X$ is \intro{$k$-tame} (or simply ``tame'' if the $k$ is either
clear form the context or unimportant)
if, for every $n \geq 1$,
for every definable subset $Y \subseteq X^n$ such that $\dim(Y) = \dim(X^n)$,
and every definable function $f: Y \to Z$, there exists $V \subseteq Y$, such
that $V$ is open in~$X^n$, large in $Y$, and $f \rest V$ is~$\Ck$.
\end{definizione}

\begin{remark}\label{rem:tame-lmin}
If $\K$ is locally o-minimal, then every definable set is $k$-tame for
every~$k$ (see~\cite{fornasiero-lomin}).
If $\K$ is not locally o-minimal, then there exists a definable subset of $\K$
that is not $0$-tame.
\end{remark}

\begin{remark}
A finitary $\Ck$-manifold is tame.
The disjoint union of finitely many tame sets (with the differential structure 
of the disjoint union) is tame.
\end{remark}

\begin{proposition}\label{prop:tame-group}
Fix $k \geq 0$.
Let $\G$ be a group definable with parameters~$A$.
Assume that $G \subseteq \K^n$ is tame and of dimension~$d$.
Let $\sigma$ be the induce topology on $G$ from~$\K^n$.
Then, there exists an $A$-definable set $V \subseteq G$ and a topology $\tau$
on~$G$, such that:
\begin{enumerate}
\item $V$ is a $\Ck$ embedded manifold;
\item $\G$ with the topology $\tau$ is a topological group;
\item $V$ is large and open in $G$ with respect both the topologies $\tau$ and
$\sigma$;
\item the topologies $\tau$ and $\sigma$ restricted to $V$ coincide;
\item $V$ is a $d$-dimensional $\Ck$ embedded manifold;
\item some $d + 1$ left (right) translates of $V$ cover~$G$:
$G = g_1 \cdot V \cup \dots \cup g_{d+1} \cdot V$;
\item
$g_1 \cdot V, \dotsc, g_{d+1} \cdot V$ 
\rom(with the obvious maps into $V$\rom) form an atlas 
inducing the topology $\tau$, and making
$G$ a finitary $\Ck$-manifold and a $\Ck$ group  
\rom(\ie, this atlas multiplication and inversion are $\Ck$-functions\rom).
\end{enumerate}
\end{proposition}
\begin{proof}
The proposition for $k = 0$ is proved in \cite{wencel}*{Theorem~3.5}.
He states the theorem under the assumption that all definable sets are
tame, but by inspecting the proof, one easily sees that it suffices that $G$
is tame.
The general proposition can be proved in the same way.
\end{proof}



Assume that $\K$ is locally o-minimal and $\G$ is a group
definable in~$\K$; then, by Remark~\ref{rem:tame-lmin}, $G$~is tame, 
and therefore we can apply Proposition~\ref{prop:tame-group} 
to $\G$  (see \cite{fornasiero-lomin}). 
This is no longer true in the case when $\K$ is not locally o-minimal.
\begin{example}
By \cite{FM05}, there exist a \dminimal expansion $\Rc$ of $\Rbar$,
and a group $\G$ definable in~$\Rc$, such
that 
\begin{enumerate}
\item
as an abstract group, $\G$ is isomorphic to 
$\pair{\Z \times \Z, +}$;
\item 
$\dim(G) = 0$;  
\item
$G$ is not $0$-tame;
\item
the conclusion of Proposition~\ref{prop:tame-group} fails for $\G$:
let $V \subset G$ be the set of isolated points of~$G$ 
\rom(in the topology induced by the ambient space\rom); 
then, $V$ is not large in $G$, and $G$ is not covered by finitely many
translates of~$V$.
\end{enumerate}
\end{example}


We can now prove a version of Pillay's result.
Remember that every definable set $X$ is the union of finitely many 
embedded manifolds (possibly, of different dimensions);
hence, after changing the topology of $X$ slightly, we can assume that $X$ is
tame; thus, we can prove the following result.
\begin{thm}\label{thm:group-topology}
Let $\pair{H, \cdot, e}$ be a definable group, with $H \subseteq \K^n$.
Then, there exist a definable group $\G$ 
and a definable continuous group isomorphism $f: G \to H$, 
such that $\G$ satisfies the conclusion of
Proposition~\ref{prop:tame-group}.
In particular, on $H$ we can put a $\Ck$ differential structure  
$\tau$ that makes $\pair{H, \cdot}$ a finitary $\Ck$-manifold and
a $\Ck$ differential group;
moreover, there exists $V \subseteq H$ definable, open in  both the topology
induced from $\K^n$ and the topology~$\tau$, 
such that $V$ is a embedded $\Ck$-manifold, and such that
$\tau$ and the $\Ck$ differential structure induced by $\K^n$
 coincide on~$V$
\rom(however, in general $V$ will not be large in~$H$, nor finitely many
translate of $V$ will suffice to cover~$H$\rom).
\end{thm}

\begin{lemma}\label{lem:unique-topology}
Let $\G$ and $\G' \coloneqq \pair{G', \cdot'}$ 
be definable groups and $k \in \N$.
\begin{enumerate}\item 
Let $\sigma$ and $\sigma'$ be topologies on $G$ and 
$G'$ respectively that make them finitary manifolds, 
and such that
all left multiplications are continuous maps.  Let $\phi: G \to G'$ be a
definable group homomorphism.  
Then, $\phi$~is continuous.
If moreover we have differential structures 
on $G$ and $G'$ that make them finitary $\Ck$-manifolds, 
such that all left multiplications are $\Ck$-functions, then
$\phi$ is also a $\Ck$-function.
\item
The  $\Ck$ differential structure  $\tau$ in Theorem~\ref{thm:group-topology} is the unique $\Ck$ differential structure on $G$ that makes it both a
a finitary manifold and a $\Ck$ differential group.
\item
Let $\phi: G \to G'$ be a definable surjective
homomorphism.
Then, $\phi$ is an open map \rom(in the group manifold topologies\rom), 
and therefore $G'$ has the quotient topology.
Moreover, if $X \subseteq G$ is clopen and definable,
then $\phi(X)$ is clopen.
\end{enumerate}
\end{lemma}
\begin{proof}
(1) Let us show the case $k = 0$.
By tameness, there exists $V \subseteq G$ open, definable, and nonempty, such
that $\phi \rest V$ is continuous.
Let $a \in G$; we have to show that $\phi$ is continuous in a neighbourhood
of~$a$.
Choose $b \in V$, and let $c \coloneqq a \cdot b^{-1}$.
Then, $c \cdot V$ is an open neighbourhood of~$b$, and, since left
multiplications by $c$ and $\phi(c^{-1})$ are continuous maps, $\phi$ is
continuous on $c \cdot V$. 
The same proof works for $k > 0$.

(2) is immediate from (1).

(3) By (1), we can assume that $\phi $ is $\Cone$.
By Sard's Lemma,  there exists $g \in G$ such that $d_g(\phi)$ has
rank equal to $\dim(G')$; thus (by the local submersion theorem) 
$\phi$ is open in a neighbourhood of~$g$.
Thus, $\phi$ is an open map.

Let $X \subseteq G$ be clopen and definable; we want to prove that $\phi(X)$
is clopen.
We know already that $\phi(X)$ is open; thus, we only need to show that it is
closed.
Since $G'$ has the quotient topology, it suffices to show that 
$Y \coloneqq \phi^{-1}(\phi(X))$ is closed in~$G$.
Let $b \in \cll(Y)$; we want to show that $b \in Y$.
Let $V \subseteq U$ be an open, definable, 
and definably connected neighbourhood of $e$ (the identity of~$\G$); 
we know that $b \cdot V \cap X \neq \emptyset$; 
but $b \cdot V$ is definably connected and $X$ is clopen and definable;
therefore, $b \in b \cdot V \subseteq X$.
\end{proof}

Thus, by the above lemma, we can talk about \emph{the} $\Ck$ differential
structure  $\tau$ that makes a definable group  $\G$ both a finitary
$\Ck$-manifold and a $\Ck$ differential group;
 we will call $\tau$ the group $\Ck$ structure  of~$\G$
(or the group manifold topology on $\G$ when $k = 0$).
When we say \eg that a definable group is definably connected, we mean in its
group manifold topology.

A similar result holds for a definable group action.
\begin{proposition}\label{prop:action}
Fix $k \in \N$.
Let $*: G \times X \to X$ be a transitive definable  group action, from a
definable group $\G$ on a definable set~$X$.
Then, there exists a $\Ck$ differential structure on $X$ such that
 $*$ is a $\Ck$-function \rom(with the group $\Ck$ structure on $\G$\rom).
\end{proposition}
\begin{proof}
See the proofs of either \cite{wencel}*{Theorem 4.6} or
\cite{PPS00}*{Theorem~2.11}.
\end{proof}

\begin{lemma}\label{lem:subgroup}
Let $\G$ be a definable group, and let $\tau$ be its group
manifold topology. 
Let $H < \G$ be a definable group.
Then, $H$ is a closed subgroup \rom(\wrt~$\tau$\rom).
Moreover, \tfae:
\begin{enumerate}
\item $H$ is clopen;
\item $H$ has nonempty interior;
\item $\dim(G/H) = 0$;
\item $\dim(G) = \dim(H)$.
\end{enumerate}
\end{lemma}
\begin{proof}
The fact that $H$ is closed and ($2 \Rightarrow 1$) are as
in~\cite{pillay87}*{Fact~2.6 and Proposition~2.7}.\\
($4 \Leftrightarrow 3$) is clear, since $\dim$ is additive.\\
($2 \Rightarrow 4$) is clear, since $\pair{G, \tau}$ is a 
finitary  manifold.\\
($4 \Rightarrow 2$) is clear, because $\dim$ is local  and
$\pair{G, \tau}$ is a finitary manifold.\\
($1 \Rightarrow 2$) is clear.\\
($2 \Rightarrow 1$) is as in~\cite{pillay87}*{Proposition~2.7}.
\end{proof}

\begin{lemma}
Fix $k \in \N$.
Let $\G$ be a definable group, with its group $\Ck$
differential structure.
Let $H < \G$ be a definable subgroup,
with its group $\Ck$ differential structure~$\tau$.
Then, $H$ is an embedded $\Ck$-submanifold of~$G$, and $\tau$ coincides with
the differential structure induced by~$G$.
\end{lemma}
Thus, when we deal with subgroups, we don't have to distinguish between the
intrinsic topology/differential structure and the induced one.
\begin{proof}
By Lemma~\ref{lem:unique-topology}, it suffices to show that $H$ is an
embedded $\Ck$-submanifold of~$G$.
Since left multiplication is a $\Ck$-function on~$G$, it suffices to show
that there exists a definable nonempty open set $V \subseteq G$, such that
$V \cap H$ is an embedded $\Ck$-submanifold of~$G$.
Let $W \subseteq G$ be an open definable subset of $G$ that is definably
$\Ck$ diffeomorphic to $\K^n$ (where $n \coloneqq \dim(G)$),
and $H' \coloneqq H \cap W$; \wloG, we can assume that $W = \K^n$.
By Fact~\ref{fact:i-regular-2},
$\reg^k(H')$ is nonempty, 
proving what we wanted. 
\end{proof}

\begin{definizione}
A definable set $X$ is \intro{definably connected} if it contains no nontrivial
\emph{definable} clopen subsets.
\end{definizione}

\begin{proposition}\label{prop:Abelian-torsion}
Let $\G \coloneqq \pair{G, +, 0}$ be a definable Abelian group, and $1 \leq n
\in \N$.
Define $\G[n] \coloneqq \set{g \in G: n g = 0}$.
Then, $\G[n]$ is a discrete subgroup of $\G$.
If moreover $G$ is definably connected, then $\G$ is divisible.
\end{proposition}
\begin{proof}
For every $1 \leq n \in \N$, let $F_n : G \to G$ be the map $x \mapsto n x$,
and $P: G \times G \to G$, $\pair{x, y} \mapsto x + y$.
Let $d \coloneqq \dim(G)$.
\begin{claim}
$d_0(F_n) = n I$, where $I$ is the $d \times d$ identity matrix in 
$M_d(\K)$, and in particular $d_0(F_n)$ is surjective.
\end{claim}
As in \cite{OPP}*{Lemma~4.3}, one sees that $d_0(P)$ 
is the matrix $\pair{I, I} \in M(2d \times d, \K)$.
Since 
\[
d_0(F_n) = d_{\pair{0,0}} P \cdot 
\begin{pmatrix}
I\\ d_0(F_{n-1})
\end{pmatrix},
\]
the claim follows by induction on~$n$.

Thus, by the local submersion theorem, $F_n$ is an open map around~$0$;
since moreover $F_n$ is a homomorphism, $F_n$ is an open map,
and in particular $F_n(G)$ is open. 
Thus, by additivity of dimension, $\dim(\G[n]) = 0$, and $\G[n]$ is discrete.

If moreover $G$ is definably connected, then, since $F_n(G)$ is an open
definable subgroup of $\G$, $F_n(G) =G$.
Since the above is true for every $n >1$, $\G$ is divisible.
\end{proof}

\begin{definizione}\label{def:dcompact}
$\G$ is \intro{definably compact} if, for every definable decreasing family
$(X_t: t \in \K)$ of closed nonempty subsets of~$G$, 
we have $\bigcap_t X_t \neq \emptyset$).
\end{definizione}

Notice that the above definition generalizes the usual one for definable subsets
of $\K^n$ with the subspace topology.

\begin{conjecture}[\cite{EO04}]
Let $\G \coloneqq \pair{G, +, 0}$ be Abelian, definably compact, definably
connected, and of dimension~$d$.
Then, for every $1 \leq n \in \N$, $G[n] \cong (\Z / n\Z)^d$.
\end{conjecture}

\begin{example}
We now show that some Lie group is not definable in any \dminimal expansion of the reals.
Fix an integer $d \geq 1$.
Let $\Qs$ be the group of nonzero rational numbers, with the usual
multiplication and the discrete topology.
Consider the action of $\Qs$ on $\R^d$ by scalar multiplication, and define
$\Hg \coloneqq \R^d \rtimes \Qs$ to be the corresponding semidirect product.
We claim that $\Hg$ is not definable in any \dminimal expansion of the real field.
In fact, assume, for a contradiction, that $\Hg$ were definable in some \dminimal structure.
Since $\Hg$ is uncountable, $\dim(\Hg) \geq 1$.
Let $\F$ be a definable Abelian subgroup of $\Hg$ of dimension at least $1$ ($\F$ exists by 
Proposition \ref{prop:Abelian-subgroup}).
It is easy to see that $\F$ is a subgroup of $\R^d \times \mset 1$; thus, its connected component $\F^0$ is a torsion-free Abelian Lie subgroup of $\R^d$ of dimension at least $1$ (unfortunately, we don't know whether $\F^0$ is definable); therefore, as a Lie group, $\F^0$ is isomorphic (not necessarily in a definable way) to $\R^e$ for some $e \geq 1$.
Let $d \coloneqq \pair{v, 1}$ be any nonzero element in $\F^0$ and $C$ be the centralizer of $d$ inside $\Hg$: notice that  $C = \set{\pair{0, q}: q \in \Qs}$.
Define, $D \coloneqq \mset{\pair{0, 1}} \cup C \cdot d = \set{\pair{q v, 1}: q \in \Raz}$.
Thus, $D$ is a definable subgroup of $\F^0$ isomorphic to $\pair{\Raz, +}$; let $R$ be the topological closure of $D$ inside $\Hg$.
By working inside the Lie group $\F^0$, we easily see that $D$ is uncountable, and therefore $\dim(D) \geq 1$, while $\dim(C) = 0$, which is impossible inside a \dminimal structure.
\end{example}

\begin{question}
Let $\Hg$ be a \emph{connected} Lie group.
Is there some \dminimal expansion $\Rc$ of the real field, such that $\Hg$ is isomorphic (as a Lie group) to a group definable in $\Rc$?
\end{question}
Notice that if $\Hg$ is compact, then the answer to be above question is yes, since we can take $\Rc$ to be the \ominimal structure $\R_{an}$; see also \S\ref{subsec:group-example}.


\subsection{Connected components}


\begin{definizione}
Let $X$ be a definable set and $a \in X$.  The \intro{definable
quasi-component} of $a$ in $X$ is the intersection of all definable clopen
subsets of $X$ containing~$a$.%
\footnote{In classical topology, given a topological space $X$ and $a \in X$,
the quasi-component of $a$ in $X$ is the intersection of all clopen subsets
of $X$ containing~$a$; it can be larger than the connected component of~$a$,
unless $X$ is locally connected.}\\  
Let $\G$ be a definable group. 
Define $\G^0$ to be the quasi-component of $e$ in the
group manifold topology.
\end{definizione}

\begin{warning}
\begin{enumerate}\item 
Let $\G$ be a definable group.  
Then, $\G^0$ is type-definable,
but, unlike in the o-minimal case, we don't know whether $\G^0$ is definable.
\item
There exists a \dminimal expansion of~$\Rbar$ that  defines
a set $X$ that is a 1-dimensional submanifold of~$\R^3$, and
such that $X$ has 2 (arc-)connected components, but it is definably connected
(see \cite{fornasiero-connected}).
Thus, unlike the o-minimal case, even for ``nice'' subsets of~$\R^n$, 
definably connected does not imply connected.
\end{enumerate}
\end{warning}

Notice that if $\G$ is a 0-dimensional definable group, then
$\G^0 = \mset{e}$.

\begin{lemma}\label{lem:G0}
Let $\G$ be a definable group, 
and $H = \G^0$.
Then, $H$ is a normal closed subgroup of~$\G$. 
If moreover $H$ is definable,
then $H$ is the smallest definable subgroup of $\G$ such that $\dim(\G/H) = 0$.
\end{lemma}
\begin{proof}
Let us show first that $H$ is a subgroup; since $x \mapsto x^{-1}$ is a
homeomorphism, it is clear that $H^{-1} = H$. Thus, we only
need to show that, given $a \in H$, $a \cdot H \subseteq H$.
Let $X \subseteq G$ be clopen and definable, such that $1 \in X$.
Notice that $a \in H \subseteq X$, and therefore
$e \in a^{-1} \cdot X$; thus, since $a^{-1} \cdot X$ 
is also clopen and definable, we have $H \subseteq a^{-1} \cdot X$,
that is, $a \cdot H \subseteq X$.
Taking the intersection of all the possible $X$, we get
$a \cdot H \subseteq H$.

The fact that $H$ is normal and closed is clear.

Moreover, since $\pair{G, \tau}$ is a finitary manifold, 
$H$~is clopen in~$G$;
thus, if $H$ is moreover definable, it must be smallest definable subgroup of $G$ such that \mbox{$\dim(\G/H) = 0$}.
\end{proof}

\begin{lemma}\label{lem:G0-type}
Let $A \subseteq \K$ and
$\G$ be a group definable over~$A$.
Then, there exists a family $\set{H_i : i \in I}$, such that
each $H_i$ is a clopen normal subgroup of $G$ definable over~$A$, 
and $\G^0 = \bigcap_{i \in I} H_i$.
\end{lemma}
\begin{proof}
Let $X \subseteq G$ be a clopen definable set, containing~$e$.
It suffices to show that there exists a clopen normal subgroup $H < G$ that is
definable over~$A$, and such that $X \subseteq H$.
\Wlog, we can assume $A = \emptyset$.
Moreover, after replacing $X$ with $X \cap X^{-1}$, we can assume that $X = X^{-1}$

Let $H(X) \coloneqq \set{g \in X: g \cdot X = X}$.
Clearly, $H(X)$ is a definable subgroup of~$G$.

\begin{claim}
$H(X) \subseteq X$.
\end{claim}
In fact, since $e \in X$, $H(X) = H(X) \cdot e \subseteq H(X) \cdot X = X$.

\begin{claim}\label{cl:H-clopen}
$H(X)$ is open (and therefore clopen).
\end{claim}
Let $U \subseteq G$ be an open, definable, and definably connected
neighbourhood of~$e$.
Let $V \subseteq U$ be an open, definable, and definably connected
neighbourhood of~$e$, such that $V^{-1} \subseteq U$.
It suffices to show that $U \subseteq H(X)$.
Let $g \in V$; we want to show that $g \cdot X = X$.
Let $x \in X$; since $X$ is clopen, we have $U \cdot x \subseteq X$; thus,
$g \cdot X \subseteq X$.
Similarly, $g^{-1} X \subseteq X$, and therefore $g \cdot X = X$.

Let $\cv$ be a finite tuple of parameters,  and $\phi(\x, \y)$ be a formula,
such that Let $\phi(G, \cv) = X$, and, for every $\cv'$,
$\phi(G, \cv')$ is a clopen subset of $G$ containing~$e$.

Let $H_0 \coloneqq \bigcap_{\cv'} H(\phi(G, \cv')$, and
$H \coloneqq \bigcap_{g \in G} g\cdot H_0 \cdot g^{-1}$.
Clearly, $H$ is a normal subgroup of~$\G$, definable without parameters,
and contained in~$X$.
Thus, it suffices to prove that $H$ is open.
Fix a parameter $\cv'$ and $g \in G$, and let
$H' \coloneqq g \cdot H(\phi(G, \cv) \cdot g^{-1}$.
By Claim~\ref{cl:H-clopen}, $g \cdot H(\phi(G, \cv)) \cdot g^{-1}$.
is clopen, and therefore contains $\G^0$.
Thus, $\G^0$ is contained in~$H$; since $\G^0$ is open in~$G$,
$H$~contains an open neighbourhood of~$e$, and therefore it is open.
\end{proof}

\begin{lemma}\label{lem:differential-hom}
Let $\G$ and 
$\G' \coloneqq \pair{G', \cdot', e'}$
be definable groups, and let $\phi$ and $\phi': G \to G'$ be definable homomorphisms.
If $d_e(\phi) = d_e(\phi')$, then $\phi \rest \G^0 = \phi' \rest \G^0$.
\end{lemma}
\begin{proof}
Same as~\cite{OPP}*{Lemma~3.2}.
The uniqueness of definable solutions to differential equations is  
Proposition~\ref{prop:differential-equation}.
\end{proof}

\begin{lemma}\label{lem:open-map-G0}
Let $\G$ and 
$\G' \coloneqq \pair{G', \cdot', e'}$
be definable groups, and let $\phi: G \to G'$ be a definable homomorphism.
Then, $\phi(\G^0) \subseteq {\G'}^0$.
If moreover $\phi$ is an open map and $\K$ is $\omega$-saturated, 
then $\phi(\G^0) = {\G'}^0$.
\end{lemma}
\begin{proof}
Let $a \in \G^0$ and $b \coloneqq \phi(a)$.
Assume, for a contradiction, that
$b \notin  {\G'}^0$; let $X \subset G'$ be clopen and definable,
such that $e' \in X$ and $b \notin X$.
Then, $Y \coloneqq \phi^{1}(X)$ is clopen and definable,
$e \in Y$, and $a \notin Y$, absurd.

Assume now that $\phi$ is open and $\K$ is $\omega$-saturated.
Let $\av \in \K^m$ be the parameters of definition of $\G$, $\G'$, and $\phi$.
Let $b \in {\G'}^0$; we want to prove that
$\phi^{-1}(b) \cap \G^0 \neq \emptyset$.
Assume not.
By Lemma~\ref{lem:G0-type}, there exists a family
$\set{H_i: i \in I}$ of clopen subgroups of~$G$, definable over~$\av$,
such that $\G^0 = \bigcap_i H_i$.
By saturation, there exists $H < G$ clopen subgroup of $G$, such that
$\phi^{-1}(b) \cap H = \emptyset$, and hence $b \notin \phi(H)$.
However,, since $\phi$ is an open map, $\phi(H)$ is a definable clopen
subgroup of $G'$, and therefore $\phi(H) \subseteq {\G'}^0$, absurd.
\end{proof}

We can now refine Proposition~\ref{prop:Abelian-torsion}
\begin{lemma}
Let $\G \coloneqq \pair{G, +,0}$ be a definable Abelian group.
Assume that $\K$ is $\omega$-saturated.
Then, $\G^0$~is divisible.
\end{lemma}
\begin{proof}
Fix $1 \leq n \in \N$, and consider the map
$\phi: G \to G$, $x \mapsto n x$.
Since $\G$ is Abelian, $\phi$ is a group homomorphism.
By Proposition~\ref{prop:Abelian-torsion}, $\phi$ is an open map.
Thus, by Lemma~\ref{lem:open-map-G0}, $\phi(\G^0) = \G^0$.
\end{proof}

We don't know if in the above lemma the assumption that $\K$ is
$\omega$-saturated is necessary.
We don't know if $\G^0$ is definable or not; however, we have the
following conjecture.

\begin{conjecture}
Let $\Rc$ be a \dminimal expansion of~$\Rbar$.
Let $X \subseteq \R^n$ be a manifold definable in~$\Rc$.
Let $Y \subseteq X$ be a clopen subset of~$X$.
Then, $\pair{\Rc, Y}$ \rom(the expansion of $\Rc$ with a predicate for~$Y$\rom)
is also \dminimal.
\end{conjecture}

\subsection{The Lie algebra of a group}

Let $\G$ be a definable group of dimension~$n$.
Following \cite{PPS00}, we can endow its tangent space 
$\gl \coloneqq T_e(G)$ with the
``usual'' Lie algebra structure, in the following way. 
For every $g \in G$, let $\chi_g : G \to G$ be the map 
$x \mapsto g \cdot x \cdot g^{-1}$.
Let $\Ad: G \to \GL_n(\K)$, $g \mapsto d_e \chi_g$ 
be the adjoint representation of~$\G$, 
and $\ad \coloneqq d_e(\Ad) : \gl \to M_n(\K)$.
Let $[\ ,\ ]$ be the Lie bracket on $T_e(G)$:
that is, $[v, w]  \coloneqq \ad(v)(w)$.
Almost everything in \cite{PPS00}[\S 2.1--2.4] goes through for \dminimal
structures, with very similar proofs (the only inconvenience is that $\G^0$
might not be definable, and hence some small changes are needed, as shown in
the proofs of the following results).

\begin{fact}\label{fact:2.31}
Let $\G$ be a definable group with Lie algebra~$\gl$.
\begin{enumerate}
\item Let $\hl$ be a linear subspace of~$\gl$.  Then, the subalgebra
$\set{v \in \gl: [v, \hl] = 0}$ is the Lie algebra of the subgroup $\set{g \in
    G: \Ad(g) \rest \hl = \id}$.
\item $\G^0$ is Abelian iff $\gl$ is Abelian \rom(that is, $[v,w] = 0$ 
$\forall v, w \in \gl$\rom).
\item If $H$ is a subgroup of $\G$, then $H^0$ is normal in $G$ iff its Lie
algebra is an ideal of~$\gl$.
  \end{enumerate}
\end{fact}

\begin{proof}
See the proofs of \cite{PPS00}*{Claims 1.31 and 1.32}.
\end{proof}

\begin{corollary}
Let $\G$ be a definable group.
If $G$ is definably connected and of dimension~$1$, then $\G$ is Abelian.
\end{corollary}
\begin{proof}
The Lie algebra of $\G$ has dimension~1, and therefore it is Abelian.
Thus, by Fact~\ref{fact:2.31}, $\G$ is Abelian.
\end{proof}

\begin{lemma}
Let $\G$ be a definable group.
Let $H$ and $L$ be definable subgroups of~$\G$.
Then, $T_e(H \cap L) = T_e H \cap T_e L $.
\end{lemma}
\begin{proof}
 $T_e(H \cap L) \subseteq T_e H  \cap T_e L $ is true for any differential manifolds $H$ and $L$.

For the opposite inclusion, we only need to show that
$\dim(L \cap H) \geq \dim(T_e  L  \cap  T_e H )$.
Consider the map $f: L \to G /H$, $l \mapsto l \cdot H$,
where on $G/H$ we put the quotient $\Cone$ structure given by Proposition \ref{prop:action}.
Let $d \coloneqq \dim(L/(L \cap H))$: notice that, for every $l \in L$,
$d = \dim L - \dim(L \cap l \cdot H)$.
By Fact \ref{fact:Sard}, there exists $l \in L$ such that $\RK(d_l f) \geq d$.
Thus, $\dim(T_l (L \cap l \cdot H)) \geq \dim(L \cap l \cdot H)$.
However, $\dim(T_l (L \cap l \cdot H)) = \dim(L \cap H)$, and
$\dim(L \cap l \cdot H) = \dim(L \cap  H)$, and we are done.
\end{proof}

Since we don't know if $\G^0$ is definable, we will use the next two lemmas.

\begin{lemma}\label{lem:center}
Let $\G$ be a definable group.
Let $H \coloneqq C_G(\G^0)$ be the centralizer of~$\G^0$, 
and $\lZ$ be the center
of $T_e(G)$, that is $\lZ \coloneqq \set{v \in T_e(G): \forall w \in T_e(G)\ [v, w] = 0}$.
Then, $H$ is definable, and $T_e(H) = \lZ$. 
\end{lemma}
\begin{proof}
Notice that $H = \set{g \in G: \Ad(g) = 0}$, and hence $H$ is
definable. By applying Fact~\ref{fact:2.31} to the subalgebra
$\mathfrak h \coloneqq T_e(G)$, we get that $\lZ = T_e (H)$.
\end{proof}

\begin{lemma}\label{lem:G0-normal}
Let $\G$ be a definable group, and let $H < \G$ be
a definable subgroup.
\begin{enumerate}
\item If $H^0$ is normal in~$\G$, then there exists a definable subgroup 
$H' < H$, such that $H^0 < H' < H$ and $H'$ is normal in $\G$.
\item If $H^0$ is Abelian and normal in~$\G$, then there exists a definable subgroup 
$H' < H$, such that $H^0 < H' < H$ and $H'$ is Abelian and normal in~$\G$.
\end{enumerate}
\end{lemma}
\begin{proof}
(1) Let $H'$ be the intersection of all $G$-conjugates of~$H$.

(2) Let $L \coloneqq C_H(H^0)$; by Lemma~\ref{lem:center}, $L$ is definable.
Let $L'$ be the center of~$L$: by assumption, $H^0 < L'$, and clearly $L'$ is
definable and Abelian.
Let $H'$  the intersection of all $G$-conjugates of~$L'$.
\end{proof}

\begin{proposition}\label{prop:Abelian-subgroup}
Let $\G$ be a definable group, and $v \in T_e(G)$.
Then, there exists a definable subgroup $H < \G$,
such that $H$ is Abelian and $v \in T_e(H)$.

In particular, if $\dim(G) \geq 1$, then there exists a definable Abelian subgroup $H < \G$, such that $\dim(H) \geq 1$.
\end{proposition}
\begin{proof}
Let $n \coloneqq \dim(G)$. 
Let $L \coloneqq \set{g \in G: \Ad(g)(v) = v}$.

By Fact~\ref{fact:2.31}, $T_e(L) = \set{w \in T_e(G): [v , w] = 0}$.

Since $[v, v] = 0$, we have $v \in T_e(L)$.
Define $M \coloneqq C_L(L^0)$ to be the centralizer of~$L^0$ inside~$L$.
\begin{claim}
$v \in T_e(M)$.
\end{claim}
By Lemma~\ref{lem:center}, $T_e(M) = \set{w \in T_e(L): [v, w] = 0}$.

It is clear that $M^0$ is an Abelian subgroup of $\G$, and that
$v \in T_e(M^0)$.
By Lemma~\ref{lem:G0-normal}, there exists $H < M$ such that $H$ is definable
and Abelian, and $M^0 < H$, and hence $v \in T_e(H)$.
\end{proof}

\begin{lemma}\label{lem:centerless}
Let $\G$ be a definable group of dimension $n$. 
Assume that $\G$ is definably connected and centerless.
Then, the adjoint map~$\Ad$ is a 
\rom(definable and $\Cinf$\rom) embedding into $\GL_n(\K)$.
\end{lemma}
\begin{proof}
The proof is in \cite{OPP}; remember that
the uniqueness of definable solutions to differential equations is  
Proposition~\ref{prop:differential-equation}.
\end{proof}

\begin{definizione}
Let $\G$ be a definable group.
We say that $\G$ is \intro{semisimple} if every  definable, normal,
Abelian subgroup is discrete.
We say that $\G$ is \intro{definably simple}
if it has no definable, normal, nontrivial subgroups.

A Lie algebra is semisimple if it has no nontrivial Abelian ideal, and it is
simple if it has no nontrivial ideal.
\end{definizione}

\begin{lemma}\label{lem:2.35}
Let $\G$ be a definable, definably connected,  semisimple group with Lie
algebra $\gl$.
Let $\hl$ be a an ideal of $\gl$.
Then, there exists a definable normal subgroup $H \trianglelefteq \G$ whose Lie algebra is~$\hl$.
\end{lemma}
\begin{proof}
As in \cite{PPS00}*{Claim~2.35}.
\end{proof}

\begin{lemma}\label{lem:2.34}
Let $\G$ be a definable and definably connected group.
\begin{enumerate}\item 
$\G$ is semisimple iff its Lie algebra is semisimple.  
\item
$\G$ is definably simple iff its Lie algebra is simple.
\end{enumerate}
\end{lemma}
\begin{proof}
(1) See the proof of \cite{PPS00}*{Theorem~2.34}, using
Lemma~\ref{lem:G0-normal} to obviate to the fact that $\Gz$ might not be
definable.

(2) See the proof of \cite{PPS00}*{Theorem~2.36}, using Lemma~\ref{lem:2.35},
and replacing everywhere ``finite'' with ``has dimension 0''.
\end{proof}

\begin{thm}
Let $\G$ be a definable group of dimension~$n$.
Assume that $\G$ is definably connected, centerless, and semisimple.
Then,
\begin{enumerate}
\item The map $\Ad: G \to \GL_n(\K)$ is an injective homomorphism, and its image is a
semi-algebraic linear group.
\item Identify $G$ with $\Ad(G) < \GL_n(\K)$.
$\G$ is the direct product of finitely many subgroups
$H_1, \dotsc H_m$,  such that
each $H_i$ is semi-algebraic, definably simple, and definably connected.
\item There exists a semi-algebraic group $\G'$ defined without parameters,
such that $\G$ is definably isomorphic to~$\G'$.
\end{enumerate}
\end{thm}
\begin{proof}
The fact that the map $\Ad$ is an injective homomorphism  is
Lemma~\ref{lem:centerless}.
The fact that the image of $\Ad$ is semi-algebraic is as in
the proof of \cite{PPS00}*{Theorem 2.37} (notice that we do have to assume that
$G$ is definably connected to conclude that $\Ad(G)$ is semi-algebraic: \eg,
if $G$ had infinitely many definably connected components, then no
definably homeomorphic copy of $G$ could be semi-algebraic.).

Thus, \wloG we can assume that $\G$ is semi-algebraic;
therefore we can work inside the structure 
$\Kb \coloneqq \pair{K, + , \cdot}$,
and use \cite{PPS00}*{Theorem~2.38} to
conclude that $\G$ is the direct product of finitely many subgroups
$H_1, \dotsc H_m$,  such that, for each $i \leq m$,
$H_i$ is semi-algebraic and definably simple in the structure~$\Kb$
(notice that we are using \cite{PPS00}*{Theorem~2.38}, not its proof).
Fix $i \leq m$.
By Lemma~\ref{lem:2.34}, the Lie algebra $T_e(H_i)$ is simple, and,
again by Lemma~\ref{lem:2.34}, $H_i$ is definably simple in the structure~$\K$.
Since $G$ is definably connected, each $H_i$ must also be definably connected.

The proof \cite{PPS02}*{Theorem~5.1} gives the third part.
\end{proof}

\subsection{Type-definable connected component}
Let $\G$ be a group definable in some $\kappa$-saturated structure, for some
large cardinal $\kappa$.
Let $\Gzz$ be the intersection of all type-definable subgroups of $\G$ of
bounded index in $\G$ (\ie, of index less than $\kappa$).
People say that ``$\Gzz$ exists'' if $\Gzz$ itself is of bounded index.
It is known that $\Gzz$ exists when $\K$ is \ominimal, and more in general
when it has NIP
(see \eg \cite{HPP08} for definitions and properties of~$\Gzz$).
In this subsection we will give an example of a \dminimal structure
such that $\Gzz$ does not exist, where
$\G$ is the group $\G \coloneqq \pair {[0,1), + \pmod 1}$.

Let $\pair{\Rbar^*, \N^*}$ be a $\kappa$-saturated elementary extension
of $\pair{\Rbar, \N}$ (of course, it is not a \dminimal structure), where
$\kappa$ is a sufficiently large cardinal.
Let $n$ be a ``non-standard natural number'', \ie $n \in \N^* \setminus \N$.
Let $P \coloneqq \set{m \in \N^*: 1 \leq m \leq n}$.
Finally, let $\K \coloneqq \pair{\Rbar^*, P}$.
Notice that $\K$ is locally \ominimal, and \emph{a fortiori} \dminimal.
We now prove that $\G^{00}$ does not exist.
Assume, for a contradiction, that $\Gzz$ exists.
Let $H < \G$ be the subgroup of infinitesimal elements; notice that $H$ is
type-definable and of bounded index in~$\G$, and therefore $\Gzz < H$.
For every $m \in P$, let $\phi_m: G \to G$ be the multiplication by~$m$,
\ie $\phi_m(g) = m g \pmod 1$.
Define $\frac H m \coloneqq \phi_m^{-1}(H)$.
Notice that $\frac H m$ is type-definable and of bounded index in~$\G$ too, 
and therefore $\Gzz < \frac H m $.
Thus, to reach a contradiction it suffices to show that
$\bigcap_{m \in P} \frac H m$ does not have bounded index in~$\G$.
Let $\lambda$ be an infinite cardinal, such that $[\G : \Gzz] < \lambda < \kappa$.
Let $\Pa{p_i: i < \lambda}$ be a sequence of elements in $\N^*$, such that,
for every $i < \lambda$,
$p_i$ is non-standard, $p_i < n$, and, for each $j < i$,
$p_i > 2 p_j p_0$.
For each $i < \lambda$, 
let $H_i \coloneqq H \cap \bigcap_{0 < j < i} \frac H {p_j}$.
It suffices to show that, for every $i < \lambda$,
$\frac H {p_i} \cap H_i$ is a proper subgroup of $\frac H {p_i}$.
By saturation, it suffices to prove the above for $i$ finite.
Thus, we have to show that, for each $i \in \N$, there exists
$g \in H \cap \frac H {p_1} \cap \dots \cap \frac H {p_i} \setminus \frac H
{p_{i+1}}$.
Let $I \coloneqq (0, \frac 1 {p_0 p_1}) \subset [0,1)$; notice that
$I \subseteq \frac H {p_1} \cap \dots \cap \frac H {p_i}$.
Notice that $\phi_{p_{i+1}}(I) = (0, p_{i+1}/{p_0 p_i}) \pmod 1 \supseteq [0,1)
\pmod 1 = G$, and therefore there exists $g \in I$ such that
$\phi_{p_{i+1}}(g) = 1/2$, and thus 
$g \in I \setminus \frac H {p_{i+1}}$, and we are done.

We conclude this subsection with a conjecture.
\begin{conjecture}
Assume that $\K$ has NIP.
Let $\G$ be a definable group.
Assume that $\G$ is definably connected and definably compact 
\rom(see Definition~\ref{def:dcompact}\rom).
Then, $\G$ has finitely satisfiable generics and satisfies compact domination
\rom(see \cite{HPP08} for the relevant definitions and properties\rom).
\end{conjecture}


\section{Definable rings}
A ring will always be associative,  
but not necessarily commutative or with~$1$;
a ring homomorphism will not necessarily send $1$ to~$1$;
a $\K$-algebra $\F$ will not necessarily contain a copy of~$\K$
(but the $1$ if $\K$ will act as the identity on~$\F$);
remember that a division $\K$ algebra is the same as a $\K$-algebra that is
also a skew field.

\begin{notation*}
In all this section, unless explicitly said otherwise, when we say that $\F$ is
a ring, we will denote by $F$ the underlying set, by $+$ the sum, by $-$ the
minus, by $\cdot$ the multiplication, by $0$ the identity of $+$, 
and by $1$ the identity of~$\cdot$, if it exists.
We define $\F^* \coloneqq F \setminus \mset 0$.
If $\F$ is definable, then $\F^0$ will be the definably connected component of
$F$ containing~$0$ (in the group topology of $\pair{F, +, 0}$).
\end{notation*}

\begin{thm}
Fix $k \in \N$.
Let $\F$ be a  ring
definable over $A \subseteq \K$, with $F \subseteq \K^n$ and $\dim(F) = d$.
Let $\tau$ be the group $\Ck$ differential structure on $\pair{F, +}$.

If $F$ is tame, then: 
\begin{enumerate}
\item $\F$ with the differential structure $\tau$ is a $\Ck$ ring;
\item we can find a definable subset $V \subseteq F$ that is large in~$F$,
open with respect to both $\tau$ and the topology  on $F$ induced by~$\K^n$, 
and an embedded $\Ck$-manifold in $\K^n$; 
\item the restriction to $V$ of $\tau$ and the $\Ck$ differentials structures  induced by $\K^n$  coincide;
\item some $d+1$ additive translates of $V$ cover $F$. 
\end{enumerate}
If moreover $\F$ is a skew field, then
the restriction of $\tau$ to $\F^*$ is the group $\Ck$ differential structure
of $\pair{\F^*, \cdot}$, and
some $d + 1$ multiplicative translates of $V \setminus \mset 0$ cover $\F^*$.

If $F$ is not tame, then there exist a definable ring 
$\F' \coloneqq \pair{F', +, \cdot, 0}$
and a definable continuous isomorphism $h: F' \to F$, such that $F'$ is tame.

In particular, there exists a definable $\Ck$ differential structure on $F$
that makes $F$ a definable $d$-dimensional $\Ck$-manifold and $\F$ a
$\Ck$ ring. 
\end{thm}
\begin{proof}
See the proof of~\cite{OPP}*{Lemma~4.1}.
Alternatively, one can adapt the proof of \cite{wencel}*{Theorem~5.1}.
\end{proof}

We would like to adapt some of the known results about rings definable in
\ominimal structures to rings definable in \dminimal structures; we will
follow the blueprints of \cites{pillay88, OPP}.

About $0$-dimensional fields, we can say almost nothing.

\begin{example}
Let $\F$ be a countable abstract ring.
Then, there exists a \dminimal expansion of $\Rbar$ that defines a
$0$-dimensional ring isomorphic to~$\F$
(see \S\ref{sec:dmin-example}\eqref{ex:fast}).
\end{example}

However, notice the following fact (which we will use later).

\begin{remark}
Let $G \leq \pair{\K, +}$ be a definable additive subgroup.  
Then, $G$ is a trivial subgroup.
\end{remark}

For higher dimensional rings instead we can say much more.
Notice that, unlike the o-minimal case, we don't have the Descending Chain
Condition for definable groups; however, we are still able to prove 
Theorem~\ref{thm:field}.

\begin{definizione}
Let $\F$ be a definable ring, and $a \in F$.
We say that $a$ is (left-)trivial if $a \cdot F = 0$;
$a$ is almost trivial if $a \cdot \F^0 = 0$.
\end{definizione}

\begin{lemma}\label{lem:d-lambda}
Let $\F$ be a definable ring  of dimension~$n \geq 1$.
Fix $1 \leq k \in \N$ and put on $\F$ the corresponding structure of 
group $\Ck$-manifold.
Define $\mu: F \to M_n(\K)$, $g \mapsto d_0(\lambda_g)$, where $\lambda_g: F
\to F$ is the left multiplication by~$g$, and $M_n(\K)$ is the ring of $n
\times n$ matrices over~$\K$.
Then, $\mu$~is a definable $\mathcal C^{k-1}$-ring homomorphism.
$\Ker \mu$ is the set of almost trivial elements of~$\F$
\rom(hence, the set of almost trivial elements of $\F$ is definable\rom). 
In particular, if
either $\F$ has no \zdvs, or $\F$ is definably connected and has no
trivial elements, then $\mu$ is an isomorphism with the image.
\end{lemma}
\begin{proof}
See \cite{OPP}*{Lemma~4.3}.
\end{proof}

\begin{lemma}\phantomsection\label{lem:linear-ring}
\begin{enumerate}
\item 
Let $G \leq \pair{\K^n, +}$ be a definable additive subgroup.
Then, $G$ is a $\K$-linear subspace of $\K^n$, and in particular it is
definably connected.
\item
Let $F \subseteq M_n(\K)$ be a definable subring \rom(not necessarily
containing~$1$\rom).  Then, $F$~is $\K$-subalgebra of $M_n(\K)$, and it is
definably connected.
\end{enumerate}
\end{lemma}
\begin{proof}
Clearly, it suffices to prove (1).
Let $c \in G \setminus \mset 0$, and define
$S_c \coloneqq \set{t \in \K: t \cdot c \in G}$.
Notice that $S_c$ is a
definable nontrivial additive subgroup of~$\K$, and therefore $S_c = \K$,
proving that $G$ is a $\K$-linear subspace of $\K^n$.
Thus, $G$~is a finite dimensional $\K$-linear space, and hence it is definably
connected.
\end{proof}

\begin{lemma}\label{lem:nonzerodivisor}
Let $\F$ be a definable $\K$-algebra.
Then, $\F$ is definably connected.
Let $a \in \F$ be a \nzdv.
Then, $\F$ has a~$1$, and $a$ has a multiplicative inverse.
\end{lemma}
\begin{proof}
By assumption, $\pair{F, +, 0}$ is a finite-dimensional $\K$-vector space,
and hence it is definably connected.
Consider the map $\lambda_a: F \to F$ (the left multiplication by~$a$).
Then, $\lambda_a$ is a $\K$-linear endomorphism of~$\pair{F, +}$; 
by assumption, $\lambda_a$
is injective, and therefore it is surjective.
Hence, there exists $u \in F$ such that $a \cdot u = a$.
Thus, for every $b \in F$, we have $a \cdot u \cdot b = a \cdot b$;
since $\lambda_a$ is injective, we have that $u \cdot b = b$ for every $B \in
F$.
Similarly, using right multiplication by $a$, we find $v \in F$ such that, for
every $b \in F$, $b \cdot v = b$.
Thus, $u = v$ is the unit of~$\F$.
Finally, since $\lambda_a$ is surjective, $a$~has a multiplicative inverse.
\end{proof}

\begin{lemma}\label{lem:connected-field}
Let $\F$ be a definable ring,
with no \zdvs, and of dimension $n \geq 1$. 
Then, $\F$~is a skew field, and, in a canonical way, a $\K$-subalgebra
of $M_n(\K)$, containing the $1$ of $M_n(\K)$, 
and it is definably connected.
\end{lemma}

\begin{proof}
Let $\mu: F \to M_n(\K)$ be the function defined in Lemma~\ref{lem:d-lambda}.
By Lemma~\ref{lem:d-lambda}, $\mu$ is a ring isomorphism; therefore, 
\wloG we can assume that $\F$ is a subring of $M_n(\K)$.
Thus, by Lemma~\ref{lem:linear-ring},
$\F$ is a $\K$-subalgebra of~$M_n(\K)$;
hence, by Lemma~\ref{lem:nonzerodivisor}, 
$\F$~is a definably connected skew field.
Moreover, by definition, $\mu(1) = 1$.
\end{proof}

Denote by $\sqrt{-1}$ 
one of the imaginary units; remember that
$\Kb$ denotes the underlying field of~$\K$.
We now state the analogue of \cite{OPP}*{Theorem~1.1}.%
\footnote{In~\cite{OPP}*{Theorem~1.1} they forgot the assumption that the ring
is infinite, which here is replaced by the assumption that is has dimension $>
0$.}

\begin{thm}\label{thm:field}
Let $\F$ be a definable  ring.
Assume that $\F$ has no \zdvs, and $\dim(F) \geq 1$.
Then $\F$ is a skew field and
\begin{enumerate}\item 
either $\dim F = 1$ and $\F$ is definably isomorphic to~$\Kb$, 
\item
or $\dim F = 2$ and $\F$ is definably isomorphic to $\Kb(\sqrt{-1})$, 
\item
or $\dim F = 4$ and $\F$
is definably isomorphic to the ring of quaternions over~$\Kb$.
\end{enumerate}
\end{thm}
\begin{proof}
By Lemma~\ref{lem:connected-field}, 
$\F$ is a finite-dimensional division $\K$-algebra (containing $\K$ in its
center). 
Conclude, as in~\cite{OPP}, by using Frobenius' Theorem.
\end{proof}




We will now study more general definable rings.
First, we will consider the definably connected ones.

First of all, notice that if $\F$ is a $\K$-subalgebra of $M_n(\K)$,
then $\F$ is definable in the language of fields 
(since it is enough to specify a $\K$-linear basis of~$\F$).

\begin{corollary}\label{cor:ring}
Let $\F$ be a definably connected definable ring of dimension $n \geq 1$.
Assume that $\F$ has no trivial elements.
Then, via the map~$\mu$, 
$\F$~is definably isomorphic to a $\K$-subalgebra of
$M_n(\K)$.
If moreover there exists $a \in \F$ that is a \nzdv, then
$\F$ contains the unit of $M_n(\K)$.
\end{corollary}
\begin{proof}
By Lemma~\ref{lem:d-lambda}, $\mu$ is an isomorphism with the image, 
and $\mu(\F)$~is a definably subring of $M_n(\K)$.
By Lemma~\ref{lem:linear-ring}, $\mu(\F)$~is a $\K$-subalgebra of $M_n(\K)$.
If $\F$ contains a \nzdv, then $\F$ contains $1$ by
Lemma~\ref{lem:nonzerodivisor} and definition of $\mu$.
\end{proof}

Thus, we have ``full'' understanding of definably connected definable rings
with no trivial elements
(\ie, each such a ring is definably isomorphic to a ``classical'' one).

\begin{lemma}\label{lem:sum}
Assume that $\K$ is $\omega$-saturated.
Let $\F$ be a definable ring, and $D \coloneqq \Ker \mu$.
Then, $F = \F^0 + D$.
\end{lemma}
\begin{proof}
Let $a \in F$.
By Lemma \ref{lem:unique-topology}, $\mu: F \to \mu(F)$ is an open map, and thus, by Lemma \ref{lem:open-map-G0} there exists $b \in \F^0$ such that $\mu(b) = \mu(a)$.
Since $a - b \in \Ker \mu$, we are done.
\end{proof}

We now give a structure theorem for definable rings with~$1$ (but not
necessarily connected).
This result is, as far as I know is new even for o-minimal structures).
\begin{thm}\label{thm:ring-1-structure}
Let $\F$ be a definable ring with~$1$.
Then, $\F^0$ is a definable subring \rom(also with~$1$, but the unit of $\F^0$
may be different from the one of $\F$\rom).
Define $D \coloneqq \Ker \mu$ and $\Dm \coloneqq \pair{D, +, \cdot, 0}$.
Then, $D$ is definable discrete subring of~$\F$ \rom(also with~$1$\rom), and
as definable rings, $\F = \Dm \rplus \F^0$:
that is, $D \cdot \F^0 = \F^0 \cdot D = 0$, $\F^0 \cap D = \mset 0$, and
$\F^0 + D = F$.

Moreover, $\mu(\F^0) = \mu(F)$; and $\mu\rest \F_0: \F^0 \to \mu(F)$ is a ring
isomorphism; thus, $\F^0$ is a $\K$-algebra.
\end{thm}

\begin{proof}
Since every $\K$-algebra is definably connected, \wloG we can assume that $\K$
is $\omega$-saturated.
\begin{claim}
$\mu(F_0) = \mu(F)$.
\end{claim}
By Lemma~\ref{lem:open-map-G0}.

In particular,  there exists $u_0 \in \F^0$ such that $\mu(u_0) = \mu(1) = 1$.
\begin{claim}
$u_0$ is a left~$1$ of $\F^0$: that is, $u_0 \cdot x = x$ for every $x \in
\F^0$.
\end{claim}
In fact, $d_0(\lambda_{u_0}) = 1 = d_0(\id_{\F_0})$;
the claim follows from Lemma~\ref{lem:differential-hom}.

By applying the same reasoning to the opposite ring of $\F$, we can conclude
that there exists $u_0' \in \F_0$ that is a right~$1$ for~$\F^0$.
However, $u_0 = u_0\ \cdot u_0' = u_0'$; and therefore $u_0$ is a $1$
of~$\F^0$.

\begin{claim}\label{cl:u-F}
$\F^0 = u_0 \cdot F$, and in particular $\F^0$ is definable.
\end{claim}
Since $\F^0$ is a bilateral ideal of $\F$, and $u_0 \in \F^0$, we have
$u_0 \cdot \F \subseteq \F^0$.
Moreover, $u_0 \cdot F \supseteq u_0 \cdot \F^0 = \F^0$.

Remember that $D = \Ker(\mu)$.

\begin{claim}
$D \cap \F^0 = (0)$, and therefore $\mu \rest \F^0$ is injective.
\end{claim}
Let $d \in D \cap \F^0$; thus, $\mu(d) = 0$, and hence $\lambda_d = 0$ on
$\F^0$; in particular, $d \cdot u_0 = 0$; but, since $d \in \F^0$,
$d \cdot u_0 = d$.

\begin{claim}
$D$ is discrete subring of $\F$ (and therefore $\dim(D) = 0$).
\end{claim}
In fact, $\Dm^0 \subseteq D \cap \F^0 = (0)$, 
and 
thus $D$ is discrete.

\begin{claim}
$D \cdot \F^0 = \F^0 \cdot D = (0)$.
\end{claim}
In fact, both $D$ and $\F^0$ are bilateral ideals of~$\F$;
thus, $D \cdot \F^0 \subseteq D \cap \F^0 = (0)$.


By Lemma~\ref{lem:sum}, $\F^0 + D = F$, and thus, as definable rings, 
$\F = \F^0 \rplus \Dm$. 
\end{proof}

Notice that the analogue of the above theorem for Lie rings is false, as the
following example show.

\begin{example}
Let $\F$ be the ring $\pair{\Z \times \R, +, \cdot}$, where
$+$ is defined component-wise, while $\cdot$ is given by
$\pair{a, b} \cdot \pair {a', b'} \coloneqq \pair{a a', a b' + b a'}$; it is
easy to verify that $\F$ is a 1-dimensional commutative ring,
with $1 = \pair{1, 0}$, and $\F^0 = \mset 0 \times \R$.
Moreover, as additive groups, $\F = \Z \rplus \R$, 
but as rings $\F \neq \Z \rplus \R$.
Notice also that $\F^0$ is a trivial ring.
The point where the proof of Theorem~\ref{thm:ring-1-structure} does not go
through is that $\mu(\pair{a, b}) = a$, and therefore
$\mu(F) = \Z$, which is not an $\R$-algebra.

Conversely, the proof of Theorem~\ref{thm:ring-1-structure} shows that
if $\F$ is a Lie ring with~$1$ and $\mu(F)$ is connected, then
$\F = \Ker \mu \rplus \F^0$ as Lie rings.
\end{example}

In general, the following construction might be useful either in finding counterexamples, or in giving structure theorems.
\begin{example}
Let $\F$ be a definable (\resp Lie) ring and $\Am$ be a definable (\resp Lie) bilateral $\F$-algebra.
Let $L \coloneqq F \times A$.
Let $+$ be the component-wise addition on $L$.
Define a product $\cdot$ in the following way:
$\pair{f, a} \cdot \pair{f', a'} \coloneqq \pair{f \cdot f', f a' + a f' + a \cdot a'}$.
Then, $\Lm \coloneqq \pair{L, +, \cdot}$ is a definable (\resp Lie) ring, and, via the identification
$\F = \F \times \mset 0$, a bilateral $\F$-algebra.
Moreover, if $\F$ has a $1$, then $\pair{1, 0}$ is the $1$ of $\Lm$.
\end{example}

We know give some partial results and conjectures for definable rings without~$1$.

\begin{proposition}
Let $\F$ be a  definable ring of dimension $n \geq 1$.
Let 
$D \coloneqq \Ker \mu$.
Assume that:
\begin{enumerate}
\item either $D$ contains no nilpotent elements;
\item or $D$ has a $1$;
\item or $D \cap \F^0 = (0)$.
\end{enumerate}
Then,
$D$ is a \rom(definable\rom) discrete subring of~$\F$.
Besides, if $\K$ is $\omega$-saturated, then 
$\F = D \rplus \F^0$ and
the map $\mu_0 \coloneqq \mu \rest \F^0$ is a ring isomorphism between 
$\F^0$ and $\mu(F)$, and hence $\F^0$ is \rom(in a canonical way\rom) 
a $\K$-algebra.
Moreover, in cases (1) and (2) $\F^0$ is definable and (3) holds.
\end{proposition}
\begin{proof}
Assume ($1$).
Let $u_1$ be the $1$ of $D$.
Let $a \in D \cap \F^0$.
Then, $a = u_1 \cdot a = 0$, and ($3$) holds.
Moreover, $\F^0 = \set{x \in F: u_1 \cdot x = 0}$, and hence it is definable.

Assume $(2)$
Let $a \in  D \cap \F^0$.
Then, $a \cdot a = 0$, and, since $D$ has no nilpotent elements, $a = 0$, and ($3$) holds.
Moreover, $\F^0 = \set{x \in F: D \cdot x = (0)}$, and hence it is definable.

Assume now ($3$).
Since $D^0 \subseteq D \cap \F^0 = (0)$, we have that $D$ is a discrete subring of $\F$.
Moreover, since both $D$ and $\F^0$ are bilateral ideals of~$\F$, $D \cdot \F^0 \subseteq D \cap \F^0 = (0)$, and similarly $\F^0 \cdot D = (0)$; thus, $\F = D \rplus \F^0$.
Since $\Ker(\mu_0) \subseteq D \cap \F^0 = (0)$, we have that $\mu_0$ is injective.
If $\K$ is $\omega$-saturated, Lemma \ref{lem:open-map-G0} implies that $\mu_0$ is also surjective.
\end{proof}



\begin{definizione}
Let $\F$ be a definable ring.
$\F$ is trivial if all element are trivial (i.e,  if $x \cdot y = 0$ for every
$x, y \in F$), and
$\F$ is almost trivial if every element is almost trivial
(\ie, if $F \cdot \F^0 = (0)$).
\end{definizione}
\begin{examples}
\begin{enumerate}\item 
If $\G$ is a definable Abelian group, then $\G$ can be made into
a trivial ring by defining $x \cdot y = 0$.  
\item
Every definable discrete ring is almost trivial.
\item
If $\F$ is a definable discrete ring, and $\G$ is a trivial ring, 
then $\F \rplus \G$ is an almost trivial ring.
\end{enumerate}
\end{examples}

\begin{lemma}\label{lem:trivial}
Let $\F$ be a definable ring.
\Tfae:
\begin{enumerate}
\item $\mu = 0$;
\item $\F$ is almost trivial;
\item the opposite of $\F$ is almost trivial
\rom(\ie, $\F^0 \cdot \F = (0)$\rom);
\item $\F^0 \cdot \F^0 = (0)$.
\end{enumerate}
\end{lemma}
\begin{proof}
($2 \Rightarrow 1$),  ($2 \Rightarrow 4$), and ($3 \Rightarrow 4$) are clear.\\
($1 \Rightarrow 2$) is Lemma~\ref{lem:d-lambda}.\\
($4 \Rightarrow 1$): by assumption, $\F^0 \subseteq \Ker \mu $; thus,
$\dim(\Ker \mu) = \dim(\F)$, and therefore $\dim(\mu(F)) = 0$;
but $\mu(F)$ is definably connected, and hence $\mu(F) = \mset 0$.\\
($3 \Rightarrow 1$): apply ($4 \Rightarrow 1$) to the opposite ring $\F^{op}$.
\end{proof}

\begin{proposition}\label{prop:ring-structure}
Let $\F$ be a definable ring of dimension $n$.
\begin{enumerate}\item 
We have a short exact sequence of definable rings
\[
0 \to \Lm \to \F \xrightarrow\mu \Am \to 0, 
\]
where 
$\Am \coloneqq \mu(\F)$ is a $\K$-subalgebra of $M_n(\K)$
and $\Lm \coloneqq \Ker \mu$ is an almost trivial ring.
\item
If $\F$ is almost trivial, then we have a short exact sequence
of definable rings $0 \to \G \to \F \to \Dm \to 0$,
where $\Dm$ is discrete and $\G$ is trivial.
\end{enumerate}
\end{proposition}
\begin{proof}
1) We only have to check that $\Lm$ is almost trivial.
However, $\Lm^0 \subseteq \F^0$, and, by Lemma~\ref{lem:d-lambda},
$\Lm \cdot \F^0 = 0$.

2) Define $G \coloneqq \set{x \in F: \forall y \in F\ x \cdot y = y \cdot
  x = 0}$.
Notice that $\F^0 \subseteq G$, and therefore $D \coloneqq \F/ G$ is discrete.
\end{proof}

\begin{open problem}
\begin{enumerate}\item 
Proposition~\ref{prop:ring-structure} shows that a definable ring is built
using a $\K$-algebra, a definable discrete ring, and a definable trivial ring.
How are these rings ``put together''?
Is $\F = \Dm \rplus \G \rplus \Am$, for some definable rings $\Dm$, $\G$,
and~$\Am$, with $\Dm$ discrete, $\G$ trivial, and $\Am$ $\K$-algebra?
\item
Let $\F$ be a definable ring \rom(not necessarily with~$1$\rom).
Is $\F^0$ definable?
Is $\F$ of the form $\Dm \rplus \F^0$,
where $\Dm$ is a definable discrete subring?
\end{enumerate}
\end{open problem}

Finally, we show that a definable discrete ring has only trivial definable connected modules.
\begin{proposition}
Let $\F$ be a definable discrete ring, and $\G$ be a definable definably
connected \rom(left\rom) $\F$-module.
Then, $\F$ acts trivially on $\G$, \ie $f g = 0$ for every $f \in F$ and
$g \in G$.
\end{proposition}
\begin{proof}
Define a ring $\Lm$ in the following way:
let $\pair{L, +} \coloneqq \pair{F, +} \times \pair{G, +}$, 
with multiplication 
$\pair{\alpha, a} \cdot \pair{\beta, b} \coloneqq \pair{\alpha \beta, \alpha b}$.
We identify $G$ with $\mset 0 \times G \subseteq L$ and $F$ with
$F \times \mset 0 \subseteq L$.
Notice that $\Lm^0 = G$, and, since $G \cdot G = 0$, Lemma~\ref{lem:trivial} 
implies that $L \cdot G = 0$, and in particular $F \cdot G = 0$, 
proving that $\G$ is a trivial $\F$-module.
\end{proof}


\bibliographystyle{alpha}	
\bibliography{tame}		

\end{document}